\documentclass[11pt,letterpaper]{article}
\pdfoutput=1 

\usepackage[T1]{fontenc}
\usepackage{datetime}
\usepackage{array,booktabs}
\usepackage[table,RGB]{xcolor}
\usepackage{arydshln}
\newcommand{\ghline}{\hdashline}
\newcommand{\cghline}{\cdashline{2-11}}
\usepackage{enumitem}

\usepackage{csvsimple}  
\usepackage{tikz}
\usetikzlibrary{arrows,calc,matrix,patterns,plotmarks,positioning}

\usepackage{hyperref}
\hypersetup{colorlinks=true,
            urlcolor=blue,
            linkcolor=black,
            citecolor=blue,
            bookmarksdepth=paragraph,
            pdftitle={Use of a Neural Network to Warmstart Column Generation for Unit Commitment},
            pdfauthor={Nagisa Sugishita, Andreas Grothey, Ken McKinnon},}

\usepackage{pgfplots, pgfplotstable}

\usepackage{amsmath, amssymb, amsthm, amsfonts}
\usepackage[font=small,labelfont=bf]{caption}
\usepackage{graphicx}

\usepackage[xindy,acronyms]{glossaries}
\glsdisablehyper  

\usepackage[
  style=authoryear,
  backend=biber,
  maxcitenames=2,
]{biblatex}

\newbibmacro*{journal}{%
  \iffieldundef{journaltitle}
    {}
    {\printtext[journaltitle]{%
       \printfield[myplain]{journaltitle}%
       \setunit{\subtitlepunct}%
       \printfield[myplain]{journalsubtitle}}}}

\DeclareFieldFormat{myplain}{#1}
\DeclareFieldFormat{titlecase}{\MakeSentenceCase*{#1}}

\usepackage[nameinlink]{cleveref}

\crefname{equation}{}{}
\crefname{section}{section}{sections}
\crefname{figure}{figure}{figures}
\crefname{table}{table}{tables}
\crefname{example}{example}{examples}
\crefname{proposition}{proposition}{propositions}
\crefname{definition}{definition}{definitions}
\crefname{theorem}{theorem}{theorems}
\crefname{remark}{remark}{remarks}

\Crefname{section}{Section}{Sections}
\Crefname{figure}{Figure}{Figures}
\Crefname{table}{Table}{Tables}
\Crefname{example}{Example}{Examples}
\Crefname{proposition}{Proposition}{Propositions}
\Crefname{definition}{Definition}{Definitions}
\Crefname{theorem}{Theorem}{Theorems}
\Crefname{remark}{Remark}{Remarks}

\theoremstyle{plain}

\theoremstyle{plain}

\theoremstyle{plain}

\theoremstyle{plain}

\theoremstyle{definition}

\theoremstyle{definition}

\theoremstyle{definition}

\theoremstyle{definition}

\theoremstyle{remark}

\numberwithin{equation}{section}

\newacronym{DW}{DW}{Dantzig-Wolfe}
\newacronym{LP}{LP}{linear programming}
\newacronym{LPR}{LPR}{linear programming relaxation}
\newacronym{MILP}{MILP}{mixed-integer linear programming}
\newacronym{MP}{MP}{master problem}
\newacronym{RMP}{RMP}{restricted master problem}
\newacronym{UC}{UC}{unit commitment}

\usepackage[author={Nagisa Sugishita},color={1.000 0.000 0.000},fontcolor={1.00 0.00 0.00},opacity=1,markup=Underline]{pdfcomment}

\newcommand{\revised}[1]{#1}
\newcommand{\revisedii}[1]{#1}
\newcommand{\revisediii}[1]{#1}

\newcommand\twoboxes[4]{%
\begin{tikzpicture}

\tikzset{block/.style={
    inner sep=.1em,
    draw=black,
    fill=white,
    minimum size=10.0cm,
    minimum height=2.8cm,
}};

\node[block] (a) {#2};

\node[block, below=0.8cm of a] (b) {#4};

\node [above left=0.1cm and 0cm of a, label={[label distance=0.1cm]right:{\bfseries\small #1}}] {};
\node [above left=0.1cm and 0cm of b, label={[label distance=0.1cm]right:{\bfseries\small #3}}] {};

\end{tikzpicture}
}

\definecolor{mysymbolcolor}{rgb}{0,0,0}

\providecommand\idxcomponent{}
\renewcommand{\idxcomponent}{{\color{mysymbolcolor} g}}

\providecommand\numcomponents{}
\renewcommand{\numcomponents}{{\color{mysymbolcolor} G}}

\providecommand\oneofidxcomponents{}
\renewcommand{\oneofidxcomponents}{%
{\color{mysymbolcolor} \idxcomponent = 1, 2, \ldots, \numcomponents}%
}

\providecommand\sumcomponents{}
\renewcommand{\sumcomponents}[1][\idxcomponent]{%
{\color{mysymbolcolor} \sum_{#1 = 1}^\numcomponents}%
}

\providecommand\randomidxcomponent{}
\renewcommand{\randomidxcomponent}{{\color{mysymbolcolor} g}}

\providecommand\dualfunction{}
\renewcommand{\dualfunction}{{\color{mysymbolcolor} q}}

\providecommand\neuralnetwork{}
\renewcommand{\neuralnetwork}{{\color{mysymbolcolor} f}}

\newcommand{\citet}[1]{\textcite{#1}}
\newcommand{\TABLE}[3]{\centering\caption{#1}{\small #2}}  

\newenvironment{APPENDICES}
  {
    \appendix
  }
  {
  }
\newcommand{\ACKNOWLEDGMENT}[1]{}

\newcommand*\samethanks[1][\value{footnote}]{\footnotemark[#1]}

\title{Use of Machine Learning Models to Warmstart\\Column Generation for Unit Commitment}
\author{Nagisa Sugishita\thanks{School of Mathematics, University of Edinburgh, James Clerk Maxwell Building, Edinburgh (UK), EH9 3FD, Email: {\tt n.sugishita@sms.ed.ac.uk}, {\tt a.grothey@ed.ac.uk}, {\tt k.mckinnon@ed.ac.uk}}, Andreas Grothey\samethanks, Ken McKinnon\samethanks}
\date{\today}

\begin{document}

\maketitle

\begin{abstract}
The unit commitment problem is an important optimization problem in the energy industry used to compute the most economical operating schedules of power plants.
Typically, this problem has to be solved repeatedly with different data but with the same problem structure.
Machine learning techniques have been applied in this context to find primal feasible solutions.
On the other hand, Dantzig-Wolfe decomposition with a column generation procedure has been shown to be successful in solving the unit commitment problem to tight tolerance.
We propose \revised{the use of machine learning models not to find primal feasible solutions directly but} to generate initial dual values for the column generation procedure.
Our numerical experiments compare \revised{machine learning based methods for warmstarting the column generation procedure} with three baselines: \revised{column pre-population,} the linear programming relaxation and coldstart.
The experiments reveal that the machine learning approaches are able to find both tight lower bounds and accurate primal feasible solutions in a shorter time compared to the baselines.
Furthermore, these approaches scale well to handle large instances.
\end{abstract}

\section{Introduction}
\label{sec_introduction}

The \gls{UC} problem is an important optimization problem in the energy industry.
Its aim is to compute the optimal operating schedules of power plants for given demand over a fixed time period.
This problem is solved by electricity generating companies on a daily basis to determine which generators are to be used.
The timings of switching the generators on and off and the amount of power dispatched have to be optimized simultaneously.
The decisions in successive time periods are coupled by ramping limits \revisedii{(the maximum rate of change in power output)} and minimum up/downtime \revisedii{(the minimum number of time periods for a generator needs to stay on/off after startup/shutdown to prevent damage)} constraints, and this gives rise to large-scale combinatorial problems.
Due to their practical importance, they have been extensively studied over the last few decades.
For a recent survey, see \citet{VanAckooijetal2018a}.

This work focuses on \gls{UC} problems that are to be solved repeatedly with different data but with the same problem structure.
This reflects practice: when a \gls{UC} problem is solved as a day-ahead planning problem, the characterisations of generators such as generation costs and ramping rates remain the same across the days, but the problems are solved with different demand forecasts each day.
This makes the problem a good candidate for the use of machine learning techniques to accelerate the solution.

\subsection{\revisedii{Literature Review}}

\revisedii{\textbf{ML for Optimization}:}
\revised{%
In recent years, the use of machine learning techniques in optimization has been studied extensively, in particular for \gls{MILP} problems.
For a survey, see \citet{BengioEtal2021}.
\revised{%
\citet{HutterEtAl2011} study automatic configuration of an \gls{MILP} solver using machine learning.
They use a local search method to find a configuration with which the \gls{MILP} solver performs well on a given set of problems.
Another popular application is the acceleration of branch and bound.
In branch and bound, the choice of branching variables has a significant impact on the overall solution time.
\citet{KhalilEtAl2016} and \citet{Gasseetal2019} train machine learning models to predict the output of Strong Branching and use the trained model as a quick surrogate of Strong Branching to select the variable to branch on.
Other applications of machine learning to branch and bound are summarised in a survey by \citet{Lodietal2017}.
}
}

\revised{An overview of applications of machine learning techniques to optimization problems specifically in the energy industry is given by \citet{YangAndWu2021}.
Various authors focus on the acceleration of the solution methods for the \gls{UC} problem.
\citet{DalalEtal2018} propose using a simple nearest neighbour method to predict the optimal cost.
In the training phase, they solve training instances to optimality and create a dataset of the optimal objective values.
Then, given a test instance, they retrieve the nearest training instance and use the corresponding optimal value as the prediction.
\citet{PinedaAndMorales2022} extend the approach to use multiple training instances and to find a feasible solution.
To solve a test instance, they retrieve a prescribed number of the nearest training instances and the corresponding optimal commitment decisions.
For each of these commitments, the binary variables in the \gls{UC} problem are fixed to the corresponding optimal solution and the remaining continuous variables are optimized.
Then, the solution with the smallest cost among the feasible solutions is adopted.
\citet{Xavieretal2020} use a modified nearest neighbour method to construct a partial solution.
Given a test instance, the average of the solutions of the nearest training instances is computed.
For each binary variable, if the corresponding average is close to 0 or 1, the variable is fixed to that value.
In this way, some of the binary variables are fixed and the resulting smaller problem is passed to an \gls{MILP} solver.
With other enhancements, the solver finds a near-optimal solution in a short time.
The focus of these studies and many of the references therein is finding good primal feasible solutions in a short time.
However, these methods do not give bounds on the suboptimality of the output.
Furthermore, they may require a large amount of time to build the training dataset: to get a single training sample, it is necessary to solve an optimization problem to optimality, and as the problem size becomes larger, the number of variables to be predicted increases, which may require larger data sets.
}

\revisedii{\textbf{Solution Methods for UC}:}
One popular traditional approach to \gls{UC} problems is to use Dantzig-Wolfe decomposition to decompose the problem by generators (e.g., see \citet{VanAckooijetal2018a}).
The reformulated problem is then solved with a column generation procedure.
This procedure can be seen as the dual of a cutting plane approach to Lagrangian relaxation, as is discussed by \citet{Briantetal2008}.
Since the \gls{UC} problem is an \gls{MILP} problem, the reformulation is not exact but a relaxation of the original problem.
However, \citet{Bertsekasetal1983} and \citet{Bard1988} reported that the integrality gap introduced by the reformulation \revised{(i.e., the difference between the optimal objective value and the lower bound provided by the relaxation)} is typically small especially if the problem size is large.
In such cases, provided that a primal heuristic finds a near-optimal feasible solution, Dantzig-Wolfe decomposition is likely to solve the \gls{UC} problem to a tight tolerance.

\revised{Recently, there has been interest in accelerating the column generation procedure using machine learning.
\citet{VaclavikEtAl2018} use a machine learning model to speed up the solver for the pricing subproblem.
They train a regression model to predict an upper bound of the objective value of the pricing subproblem, and this upper bound is passed to the optimization solver.
\citet{ShenEtAl2022} study the column generation procedure applied to a graph colouring problem.
They use a machine learning model to generate a near-optimal solution to the pricing subproblem.
When their method fails to find a solution to the pricing subproblem that has negative reduced cost, an optimization solver is used to solve it exactly.
\citet{Morabit2021} use a machine learning model to select columns from equally-promising \revisedii{ones obtained by solving} the pricing subproblem and add them to the \gls{RMP}.
Their approach is especially of value when primal degeneracy is present in the problem.
}

\revisedii{\textbf{Warmstarted Column Generation}:}
In the aforementioned decomposition-based approaches, the dual values play a key role, and it is expected that using appropriate dual values as an initial point will speed up the solution method.
One approach to generating such dual values is to solve an approximation of the original problem and obtain its optimal dual values.
\citet{Borghettietal2002} and \citet{Schulzeetal2017} relax the integrality constraints and obtain a continuous relaxation.
However, the relaxation has a similar number of variables and constraints as the original problem, and solving it even without the integrality constraints takes significant computational time.
\citet{Takritietal1996} drop further constraints such as minimum up/down time constraints and minimum power output constraints.

\revised{A different approach to warmstarting the column generation procedure is to pre-populate columns in the \gls{RMP} as is discussed by \citet{VanHoaiEtAl2005}.
In this approach, ``useful'' columns are added to the \gls{RMP} first, which requires being able to generate and identify useful columns in advance.
In some applications, such columns may be generated using domain-specific knowledge.
Alternatively, if a family of \gls{MILP} problems is being solved, the columns generated for previously solved similar instances may be added.
}

\revised{%
As we will see in our numerical experiments, when the column generation procedure is applied to the \gls{UC} problem, these initialisation methods (e.g., solution of the \gls{LPR}) can account for a significant part of the total computational time and the quality of the initialisation has a big effect on the speed of convergence of the column generation procedure.
In this paper, to speed up the initialisation time,} we propose to warmstart the column generation procedure using machine learning techniques.
Namely, we use a machine learning model to generate initial dual values.
A model is first trained to output dual values which yield a tight dual lower bound.
After the training, when solving a new instance, the machine learning model is used with the problem parameter as input to generate dual values.
The generated dual values are then used to warmstart the column generation procedure and the column generation procedure allows us to further tighten the lower bound and, with the aid of a primal heuristic, obtain feasible solutions.
With this approach, we can exploit the strength of the machine learning techniques while maintaining the desirable property of the column generation procedure (with suitable primal heuristics), such as the capability to provide high-quality solutions with very tight, provable lower bounds.
\revised{%
Furthermore, the dimension of the dual variables is much smaller than that of the primal variables and does not depend on the number of generators but only on the number of time periods.
Therefore we may expect that learning the dual values will be easier than learning primal solutions directly.
}

\revised{One approach is to train a machine learning model to predict the optimal dual values from the problem parameter values using supervised learning.
In the training phase, a solution method such as the column generation procedure is run on training instances to create a dataset of optimal dual values.
Then, a regression model is trained using the dataset.
}%
\revisedii{
This is closely related to the approach studied by \citet{PinedaAndMorales2022}: they train machine learning models to predict the optimal primal solution using a dataset of optimal primal solutions, whereas we train the model to predict the optimal dual values using a dataset of optimal dual values.
}

\revised{In an alternative approach proposed by \citet{Nairetal2018}, a neural network is trained to maximise the dual lower bound directly, without relying on a pre-built training dataset.
In their study, a dual decomposition is applied to a parametrized two-stage stochastic programming problem.
Using the decomposable structure of the dual lower bound by scenarios, an efficient stochastic gradient-based method is devised to train the neural network.
This only requires the solution of a single scenario subproblem in each iteration and does not require a dataset of optimal dual values.
}

\revisedii{
Another possible approach is based on a surrogate model.
In the approach we tested, the surrogate is a regression model trained to predict the value of the dual function given the problem parameter and dual values.
When a test instance is given, the problem parameter is fixed to that of the test instance, and the dual values are varied to maximise the regression model.
The dual values found in this way are used to initialise the column generation procedure.
\revisediii{One drawback of this approach is the computational time.
To obtain the dual values we need to solve an optimization problem.
Although the surrogate model is cheaper to optimize than the original problem, we found it still takes a longer time than the other approaches.
Given this limitation,} we do not consider this approach in this paper.
\revisediii{For detailed discussions on this topic,} see the survey paper by \citet{Jones2001}.
}

\revisedii{\textbf{Contributions}:}
\revised{The main contributions of this paper are as follows.
Firstly, we demonstrate the use of machine learning to predict dual values that can be used to initialise the column generation procedure.
By combining the machine learning techniques and the column generation procedure, we can exploit the strength of the two; namely, fast evaluation of machine learning techniques and high accuracy of the column generation procedure with provable suboptimality.
This provable suboptimality property is missing in many earlier applications of machine learning to the \gls{UC} problem.
Secondly, we provide comprehensive numerical experiments to compare the performance of the proposed approaches on large-scale instances.
The performance is measured both in terms of the tightness of the initial lower bounds and in terms of the solution time required by the warmstarted column generation procedure to find an accurate primal solution of proven suboptimality.
}

The rest of the paper is structured as follows.
Section \ref{sec_dantzig_wolfe_decomposition} briefly reviews Dantzig-Wolfe decomposition and the column generation procedure.
Section \ref{sec_warmstarting} presents methods based on machine learning models to generate initial values for the column generation procedure.
In Section \ref{sec_neumerical_experiments} the proposed approach is applied to large-scale \gls{UC} problems.
Finally, in Section \ref{sec_conclusion} conclusions and further extensions of this work are presented.

\section{Dantzig-Wolfe Decomposition and Column Generation}
\label{sec_dantzig_wolfe_decomposition}

\glsreset{RMP}
\glsreset{LPR}

\revised{As noted in Section \ref{sec_introduction}, Dantzig-Wolfe decomposition is known to be effective for the \gls{UC} problem.}
In this section, we briefly review Dantzig-Wolfe decomposition and the column generation procedure.
For further background, see \citet{VanderbeckAndSavelsbergh2006}.

Consider the following family of \gls{MILP} problems parametrized by $\omega$:
\begin{align}
z(\omega) = \min_{x_1, \ldots, x_G} \ & \sumcomponents c_\idxcomponent^T x_\idxcomponent
\label{problem_eq_mip} \\
\text{s.t.} \ &
\sumcomponents A_\idxcomponent x_\idxcomponent = a(\omega), \notag \\
& x_\idxcomponent \in X_\idxcomponent,
\quad \oneofidxcomponents, \notag
\end{align}
where \revised{$G$ is the number of subproblems,}
$x_1, x_2, \ldots, x_\numcomponents$ are vectors of decision variables and
\[
X_\idxcomponent = \{
  x_\idxcomponent \in \{0, 1\}^{n} \times \mathbb{R}^{m} \mid D_g x_\idxcomponent \le d_g
\},
\quad \oneofidxcomponents.
\]
\revised{%
Here, \revisedii{$n$ and $m$ are the number of integer and continuous variables respectively in $x_g$, and \revisediii{for each $g$ and $\omega$} $c_g \in \mathbb{R}^{n + m}$, $A_g \in \mathbb{R}^{k \times (n + m)}$, $a(\omega) \in \mathbb{R}^k, D_g \in \mathbb{R}^{l \times (n + m)}, d_g \in \mathbb{R}^{l}$, \revisediii{where $k$ and $l$ denote the number in each of the corresponding constraints}.}
}
We assume that $X_\idxcomponent$ is non-empty and bounded for $\oneofidxcomponents$, and that the problem \eqref{problem_eq_mip} has a feasible solution for every $\omega$.
To reduce clutter, in what follows we drop the dependence on $\omega$ except where this might cause confusion.
We also write $x = (x_1, x_2, \ldots, x_G)$ and $X = X_1 \times X_2 \times \cdots \times X_G$.
\revised{%
In the \gls{UC} problem, $G$ is the number of generators,  $x_g \in X_g$ corresponds to a feasible operational plan of generator $g$, $\omega$ is the vector of demands and the first constraint in \eqref{problem_eq_mip} represents the system-wide constraints (i.e., the load balance and spinning reserve constraint).
The constraint right-hand-side $a(\omega)$ has entries for the demand and spinning reserve.
The complete formulation of the \gls{UC} problem is given in Appendix \ref{sec_appendix_problem_formulation}.
}

In Dantzig-Wolfe decomposition, we consider a relaxation of \eqref{problem_eq_mip}, referred to as the \gls{MP}, in which $X_\idxcomponent$ is replaced by $\mathrm{conv}({X_\idxcomponent})$, the convex hull of $X_\idxcomponent$, for every $\idxcomponent$.
Given the boundedness assumption on $X_\idxcomponent$ the \gls{MP} can be written in terms of the extreme points $\{x_{\idxcomponent i} \mid i \in I_\idxcomponent\}$ of $X_\idxcomponent$ as
\begin{align}
\min_{p} \ & \sumcomponents \sum_{i \in I_\idxcomponent} c_\idxcomponent^T x_{\idxcomponent i} p_{\idxcomponent i}
\label{eq_mp} \\
\text{s.t.} \ &
\sumcomponents \sum_{i \in I_\idxcomponent} A_\idxcomponent x_{\idxcomponent i} p_{\idxcomponent i} = a,
\label{eq_mp_constraint_one} \\
& \sum_{i \in I_\idxcomponent} p_{\idxcomponent i} = 1, \quad \oneofidxcomponents,
\label{eq_mp_constraint_two} \\
& p_{\idxcomponent i} \ge 0, \quad \oneofidxcomponents, i \in I_\idxcomponent. \notag
\end{align}
\revised{%
This is a \gls{LP} problem with decision variables $p_{\idxcomponent i}$ ($\oneofidxcomponents, i \in I_\idxcomponent$).
}

Since the \gls{MP} is usually too large to formulate and solve explicitly, a column generation procedure is used.
The \gls{RMP} is defined by replacing each $I_\idxcomponent$ in the \gls{MP} with a subset $\hat{I}_\idxcomponent \subset I_\idxcomponent$.
We assume that the \gls{RMP} is feasible (the sets $\hat{I}_\idxcomponent$ must have suitable columns to satisfy the requirement of constraint \eqref{eq_mp_constraint_one}) and define $y$ and $\sigma_\idxcomponent$ for $\oneofidxcomponents$ to be the optimal dual values to the \gls{RMP} corresponding to the restricted version of constraints \eqref{eq_mp_constraint_one} and \eqref{eq_mp_constraint_two}, respectively.
To find columns to be added to the \gls{RMP}, the following subproblems, known as the pricing subproblems, are solved:
\begin{equation}
r_\idxcomponent(y) = \min_{
  x_\idxcomponent} \{(c_\idxcomponent^T - y^T A_\idxcomponent) x_\idxcomponent \mid x_\idxcomponent \in X_\idxcomponent
\}, \quad \oneofidxcomponents.
\label{problem_eq_pricing_subproblem}
\end{equation}
\revised{In the \gls{UC} problem, this pricing subproblem is a scheduling problem of a single generator: given a \revisedii{reduced} cost $c_g^T - y^T A_g$, it finds the cheapest operational plan for the corresponding generator.}
If $r_\idxcomponent(y) \ge \sigma_\idxcomponent$ for all $\idxcomponent$, the \gls{RMP} already includes all relevant columns and has found the optimal solution to the \gls{MP}.
Otherwise, the solutions to the pricing subproblems for which $r_\idxcomponent(y) < \sigma_\idxcomponent$ are added to the \gls{RMP} and the above process is repeated.

\revised{%
For each $y$, a lower bound on the optimal objective value of \eqref{problem_eq_mip} can be obtained using duality.
Let us denote the Lagrangian of \eqref{problem_eq_mip} by
\begin{equation*}
L(x, y)
= \sum_{g = 1}^G c_g^T x_g - y^T \left(\sum_{g = 1}^G A_g x_g - a \right)
= a^T y + \sum_{g = 1}^G \left(c_g^T - y A_g \right) x_g.
\end{equation*}
Since the minimisation \revisedii{of $L(x, y)$ over $x \in X$ for any $y$} is a relaxation of \eqref{problem_eq_mip}, the the optimal objective value $z$ of \eqref{problem_eq_mip} satisfies
\begin{equation}
z \ge \min_{x \in X} L(x, y) =: q(y),
\label{problem_eq_definition_of_q}
\end{equation}
for any $y$.
From \eqref{problem_eq_pricing_subproblem} and the definition of the Lagrangian, it follows that
\begin{equation}
\dualfunction(y) = a^T y + \sumcomponents r_\idxcomponent(y).
\label{problem_eq_dw_lower_bound}
\end{equation}
We refer to this value as the {\it dual lower bound}.
}

\revised{
It is known that the simple column generation approach discussed above suffers from instability.
In the first few iterations, due to the poor \gls{RMP} model, the column generation procedure tends to output irrelevant dual values.
\citet{Vanderbeck2005} refers to this behaviour as the ``heading-in effect''.
This is closely related to a well-known instability issue of the cutting-plane algorithm.
See for example the discussion in Chapter XII and XV in the book by \citet{HiriartUrrutyAndLemarechal1993b}.
To tackle the issue, it is necessary to deploy some mechanism to stabilise the dual values, such as quadratic regularisation of the dual values or the box step method as described by \citet{Briantetal2008}.
In this work, quadratic regularisation of the dual values is used.
\revisedii{We note that regularisation is also used to cope with degeneracy when the RMP is degenerate.
However, in our case, we did not observe any degeneracy of the \gls{RMP} when solving the \gls{UC} problem.
The purpose of adding the regularisation in our context is to mitigate the ``heading-in effect''.}
}

Dualizing the \gls{RMP} and adding quadratic regularisation on the dual values gives
\begin{align}
\max_{y, \sigma} \ & a^T y + \sumcomponents \sigma_\idxcomponent - \frac{\mu}{2} \| y - \bar{y} \|_2^2
\label{eq_regularised_rmp} \\
\text{s.t.} \ &
\sigma_\idxcomponent \le (c_\idxcomponent^T - y^T A_\idxcomponent) x_{\idxcomponent i} \quad \oneofidxcomponents, i \in \hat{I}_\idxcomponent, \notag \\
& y, \sigma: \mathrm{free}, \notag
\end{align}
where $\bar{y}$ is a regularisation centre and $\mu$ is a parameter used to adjust the strength of the regularisation.
\revised{%
We refer to \eqref{eq_regularised_rmp} as the regularised \gls{RMP}.
}
In the regularised column generation procedure, the regularised \gls{RMP} is used in place of the \gls{RMP}.
The optimal solution to \eqref{eq_regularised_rmp} is computed and the pricing subproblems are solved based on this solution.
The regularisation centre $\bar{y}$ is updated to the current dual values $y$ whenever the lower bound \eqref{problem_eq_dw_lower_bound} has improved.
In the following, we normally omit the word ``regularised'' when we refer to the regularised \gls{RMP} or the regularised column generation procedure but add the word ``unregularised'' when referring to the one without regularisation (i.e., $\mu = 0$).

\subsection{Warmstarting the Column Generation Procedure}
\label{sec_dantzig_wolfe_decomposition_subsec_warmstarting_the_column_generation_procedure}

In the above algorithm, the dual values $y$ are updated iteratively to provide tighter lower bounds.
It can be shown that the lower bound $\dualfunction(y)$ is continuous (in fact concave) in $y$.
It follows that in the case where near-optimal dual values are used as the initial point of the algorithm, the lower bound in the first iteration will be close to the tightest.
Assuming that the lower bound is sufficiently tight, the algorithm would terminate as soon as a primal feasible solution with sufficiently small suboptimality was found by primal heuristics.
It is therefore expected that using near-optimal dual values for the initial dual values will help the algorithm to terminate in a shorter time.

It is worth noting that the regularisation of the dual values as discussed above is crucial in this context.
Without regularisation, even if near-optimal dual values are used as the initial point, the algorithm is likely to be quite unstable, yielding dual values with poor lower bounds in the following iterations.
See \citet{Briantetal2008} for further discussion.

As discussed in Section \ref{sec_introduction}, one approach to finding good dual values is to solve an approximation of \eqref{problem_eq_mip}.
\revised{%
Let us denote the linear relaxation of $X_g$ by $\bar{X}_g$:
\begin{equation*}
\bar{X}_g = \{
  x_g \in [0, 1]^{n} \times \mathbb{R}^{m} \mid D_g x_g \le d_g
\}, \quad g = 1, 2, \ldots, G.
\end{equation*}
The \gls{LPR} is obtained by replacing $X_g$ with $\bar{X}_g$ in \eqref{problem_eq_mip}.
}
\citet{Borghettietal2002} and \citet{Schulzeetal2017} solve the \gls{LPR} to optimality and use the optimal dual values to the \gls{LPR} as the initial dual values for the column generation procedure.
This approach does not require any training as is required in a machine learning model.
On the other hand, solving the \gls{LPR} takes a non-trivial amount of time.

\revised{Another approach is to initialise an unregularised \gls{RMP} with ``useful'' columns and use the optimal dual values to the unregularised \gls{RMP} as the regularisation centre of the following iteration.
To this end, we need to be able to generate useful columns before running the column generation procedure.
For example, if we solve similar instances sequentially, we may use the columns generated for the previous instances.

}

\subsection{Initial Dual Lower Bounds in the Column Generation Procedure}
\label{sec_dantzig_wolfe_decomposition_subsec_initial_dual_lower_bounds_in_cg}

\revised{%
We note that the \gls{LPR} is a relaxation of \eqref{problem_eq_mip} so it gives a lower bound on the optimal objective value.
On the other hand, any dual values give a lower bound on the optimal objective value of \eqref{problem_eq_mip} by \eqref{problem_eq_dw_lower_bound}.
In particular, the optimal dual values to the \gls{LPR} give such dual values.
In the first iteration, the column generation procedure with the \gls{LPR} initialisation evaluates this lower bound.
It is of interest to compare these two lower bounds.
Let $\bar{x}^*$ and $\bar{y}^*$ be any optimal primal and dual values to the \gls{LPR}, respectively, and let $\bar{z}^*$ be the optimal primal and dual objective value of the \gls{LPR}.
Using \eqref{problem_eq_definition_of_q} and the fact that $\bar{X} \supset X$, we have
\begin{equation*}
q(\bar{y}^*)
= \min_{x \in X} L(x, \bar{y}^*)
\ge \min_{x \in \bar{X}} L(x, \bar{y}^*).
\end{equation*}
On the other hand, since $\bar{y}^*$ is the optimal dual values to the \gls{LPR}, strong duality (e.g., see \citet{Bertsekas2009}) gives
\begin{equation*}
\min_{x \in \bar{X}} L(x, \bar{y}^*)
= L(\bar{x}^*, \bar{y}^*)
= \bar{z}^*.
\end{equation*}
Thus, the dual lower bound $q(\bar{y}^*)$ computed using the optimal dual values to the \gls{LPR} is at least as tight as (and probably tighter than) the optimal objective value of the \gls{LPR}, $\bar{z}^*$, i.e.,
\begin{equation*}
q(\bar{y}^*) \ge \bar{z}^*.
\end{equation*}
}

\section{Machine Learning Methods to Compute Initial Dual Values}
\label{sec_warmstarting}

\revised{%
In Section \ref{sec_dantzig_wolfe_decomposition_subsec_warmstarting_the_column_generation_procedure}, we briefly discussed existing approaches to warmstarting the column generation procedure.
In this section, we consider approaches to training a machine learning model to generate initial dual values for the column generation procedure.
In a practical situation, the goal is to solve the \gls{UC} problem as quickly as possible.
As we will observe in our numerical experiments, there is a strong connection between solution time and the tightness of the initial dual lower bound.
Thus, we will train a machine learning model to output dual values that yield a tight dual lower bound \eqref{problem_eq_dw_lower_bound}, using this as a surrogate measure of solution time.
Below, we consider two approaches to achieve this goal.
}

\subsection{Machine Learning Model based on Single-Sampling Training}
\label{subsec_machine_learning_model_based_on_single_sampling}

\revised{%
A simple approach to training a machine learning model is to build a dataset of optimal dual values of training instances and train a regression model to predict the optimal dual values from the problem parameter.
To build a dataset\revisedii{,} a set of training instances must be solved for example with the column generation procedure.
This approach based on a dataset of optimal solutions is used to predict optimal solutions of optimal power flow problems by \citet{GuhaEtAl2019, ZamzamAndBaker2020, OwerkoEtAl2020}.
One can use any prediction model to this end.
In our numerical experiments\revisedii{,} we will use the alternatives of a neural network model, a random forest model and a nearest neighbour model.
\revised{%
We refer to this approach in the remainder of the paper as {\it single-sampling training} since this only involves the sampling of the problem parameter $\omega$.
}
}

\subsection{Machine Learning Model based on Double-Sampling}
\label{subsec_machine_learning_model_based_on_double_sampling}

\revised{%
A potential drawback of the single-sampling training is the large amount of time required to solve enough problems to build a sufficiently large dataset.
This is especially problematic when the problem is large.
An alternative approach is to train a machine learning model to maximise the expected dual lower bound directly.
This approach can exploit the decomposable structure of the dual function.
}
This approach was introduced by \citet{Nairetal2018} and used to solve a parametrized two-stage stochastic programming problem.

A neural network model can be seen as a function $\neuralnetwork(\omega, \theta) = y$ which maps the problem parameter $\omega$ and the model parameter $\theta$ to a dual value $y$.
We aim to learn values of the model parameter $\theta$ so that given $\omega$ the model outputs a dual value $y$ for which the dual lower bound $\dualfunction(y)$ is tight.
In this section, we explicitly show the dependency of the lower bound on $\omega$ as $\dualfunction(\omega, y)$, which was suppressed in \eqref{problem_eq_dw_lower_bound}.
Assume that the distribution of $\omega$ is given (e.g., the empirical distribution based on historical data).
Our goal is to maximise the expected lower bound
\begin{align*}
p(\theta) = \mathrm{E}_\omega[\dualfunction(\omega, \neuralnetwork(\omega, \theta))].
\end{align*}
From \eqref{problem_eq_dw_lower_bound}, it follows that
\begin{align*}
\dualfunction(\omega, y)
= a(\omega)^T y + \sumcomponents r_\idxcomponent(y) = \frac{1}{|\numcomponents|} \sumcomponents (a(\omega)^T y + |\numcomponents| r_\idxcomponent(y)).
\end{align*}
We can interpret the final term as the expectation $\mathrm{E}_{\randomidxcomponent}[\widetilde{\dualfunction}(\omega, y, \randomidxcomponent)]$
where
$g$ is uniformly sampled from $\{1, 2, \ldots, G\}$ and
\begin{align*}
\widetilde{\dualfunction}(\omega, y, \randomidxcomponent) = a(\omega)^T y + |\numcomponents| r_\randomidxcomponent(y).
\end{align*}
If follows that
\begin{align*}
p(\theta) = \mathrm{E}_{\omega, \randomidxcomponent}[\widetilde{\dualfunction}(\omega, \neuralnetwork(\omega, \theta), \randomidxcomponent)].
\end{align*}
$p(\theta)$ can be seen as an expectation over both $\omega$ and $g$.
We use the stochastic gradient ascent method, the standard approach used for training a neural network.
That is, we sample $\omega$ and $\randomidxcomponent$, compute the gradient of $\widetilde{\dualfunction}$ with respect to $\theta$ and make a single gradient ascent step.
We then resample $\omega$ and $\randomidxcomponent$ and repeat the process.

The gradient of $\widetilde{\dualfunction}$ with respect to $\theta$ can be computed as follows: fix problem parameter $\omega$ and subproblem index $\randomidxcomponent$ and compute the dual values $y = \neuralnetwork(\omega, \theta)$ and the component $\widetilde{\dualfunction}(\omega, y, \randomidxcomponent)$ of the dual lower bound corresponding to subproblem $\randomidxcomponent$.
Suppose that the model output $\neuralnetwork(\omega, \theta)$ is differentiable with respect to $\theta$ and the optimal value $r_\randomidxcomponent(y)$ of the pricing subproblem \eqref{problem_eq_pricing_subproblem} is differentiable with respect to $y$ \revised{(which is the case when the optimal solution to the pricing problem is unique)}.
Then the gradient of $\widetilde{\dualfunction}(\omega, y, \randomidxcomponent)$ with respect to $y$ is given by
\begin{equation}
\frac{\partial \widetilde{\dualfunction}}{\partial y} = a - |\numcomponents| A_\randomidxcomponent x^*_\randomidxcomponent
\label{eq_gradient_of_dual_function_component}
\end{equation}
where $x^*_\randomidxcomponent$ is the solution to the pricing subproblem $\randomidxcomponent$ \eqref{problem_eq_pricing_subproblem}.
Using the chain rule, we obtain
\begin{equation}
\frac{\partial \widetilde{\dualfunction}}{\partial \theta} = J \frac{\partial \widetilde{\dualfunction}}{\partial y},
\label{eq_sampled_gradient_dl_dtheta}
\end{equation}
where $J$ is the Jacobian matrix of the neural network output $y = \neuralnetwork(\omega, \theta)$ with respect to $\theta$, which is given by automatic differentiation.
\revised{%
In this paper, we refer to this approach as {\it double-sampling training} since this involves sampling both the problem parameter $\omega$ and the subproblem $g$.
}

\revised{%
The procedures of the single and double-sampling training are listed in Figure \ref{warmstart_fig_training_procedures}.
It is important to note that each step of the double-sampling training involves evaluation of \eqref{eq_sampled_gradient_dl_dtheta}, which requires the solution of only a single pricing subproblem, which is typically substantially smaller than the original problem.
On the other hand, the single-sampling training requires solving the training \gls{UC} instances to optimality.%
}
\revisedii{%
Furthermore, the single-sampling training uses the mean squared error to train a model, that is, the model is encouraged to output the dual values close to the optimal dual values.
However, the dual values close to the optimal values in terms of the Euclidean distance do not necessarily lead to a tight dual lower bound.
In contrast, the double-sampling training uses the dual lower bound to train the model, which is the metric we are directly interested in.
}

\revisedii{
One advantage of the single-sampling training is its flexibility.
The training is based on the dataset of the optimal dual values.
Once the dataset is built, one can fit models with different architectures.
}

\begin{figure}[htbp]
\centering  
\twoboxes{Single-Sampling Training}{%
    \begin{minipage}{15.5cm}
    \small
    \begin{itemize}[noitemsep]
    \item Solve as many training instances of \eqref{problem_eq_mip} with different $\omega$ as possible within the training budget.
    \begin{itemize}[noitemsep]
        \item Save each problem parameter $\omega$ and the corresponding near-optimal dual values.
    \end{itemize}
    \item Train a regression model.
    \begin{itemize}[noitemsep]
        \item Train to predict the optimal dual values from given $\omega$.
    \end{itemize}
    \end{itemize}
    \end{minipage}
}{Double-Sampling Training}{%
    \begin{minipage}{15.5cm}
    \small
    \begin{itemize}[noitemsep]
    \item Repeat as often as possible within the training budget:
    \begin{itemize}[noitemsep]
        \item Sample problem parameter $\omega$ and generator $g$.
        \item Compute the gradient of $\widetilde{q}(\omega, f(\omega, \theta), g)$ with respect to $\theta$ using \eqref{eq_gradient_of_dual_function_component} and \eqref{eq_sampled_gradient_dl_dtheta}.
        \item Do a stochastic gradient step to improve $p(\theta)$.
    \end{itemize}
    \end{itemize}
    \end{minipage}
}
\caption{\revised{Diagrams to show the training procedures. The upper one corresponds to the single-sampling training while the lower one is to the double-sampling training. \label{warmstart_fig_training_procedures}}}
\end{figure}

Once a model is trained with the single or double-sampling training, it can be used to compute initial dual values for the regularised column generation procedure, and it is expected that unlike solving an approximation such as the \gls{LPR} the computation will be quick.
Furthermore, we later observe experimentally that the trained model produces high-quality initial dual values, in terms of both the tightness of the resulting dual lower bound and the solution time of the column generation procedure when warmstarted from it.

\section{Numerical Experiments}
\label{sec_neumerical_experiments}

In this section, the performance of the dual initialisation methods based on machine learning models as well as the benchmark initialisation methods is evaluated on \gls{UC} problems of various sizes.\footnote{The source code is available at: \href{https://github.com/nsugishita/ml_to_warmstart_cg}{https://github.com/nsugishita/ml\_to\_warmstart\_cg}}
All methods are implemented in Python.
IBM ILOG CPLEX 20.1.0\footnote{\href{https://www.ibm.com/products/ilog-cplex-optimization-studio}{https://www.ibm.com/products/ilog-cplex-optimization-studio}} is used as the optimization solver (the barrier method for the \gls{RMP} and the branch and bound method for the pricing subproblems), and PyTorch and scikit-learn are used to implement the neural network models and the random forest models, respectively.  
The experiments are performed on a workstation with a 16-core Intel\textsuperscript{\tiny\textregistered} Xeon\textsuperscript{\tiny\textregistered} E5-2670 and 126 GB of RAM.

\subsection{Problem}

In the experiments, we consider a setup in which \gls{UC} problems are solved repeatedly with a fixed set of generators but with different demand forecasts.
To assess the scalability, we consider 3 different problem sizes, i.e., problems with 200, 600 and 1,000 generators.
In all cases, the length of the planning horizon is 48 hours with a time resolution of 1 hour.
The generator data is based on \citet{Borghettietal2002}.
Since their sets of generators contain 200 generators at most, we combine multiple sets to create larger ones.
For example, to create a \gls{UC} instance with 1,000 generators, we combine 5 distinct 200-generator sets.
Each generator is unique \revised{and distinct across the sets} so combining these sets does not introduce symmetry.
The demand data is based on the historical demand data in the UK published by National Grid ESO.\footnote{\href{https://www.nationalgrideso.com/}{https://www.nationalgrideso.com/}}
A detailed description of the problem formulation and implementation details of the initialisation methods are given in Appendix \ref{sec_appendix_problem_formulation}.

\subsection{Initialisation Methods}

The initialisation methods used in the experiments are described below.

\subsubsection{Benchmark initialisation Methods.}

The first three methods are not based on machine learning models but are evaluated as benchmarks.

\paragraph{Coldstart:}
as a naive baseline, column generation is run from initial dual values $y = 0$.
This method does not require any training.

\paragraph{LPR:}
this is the method based on the \gls{LPR} described in Section \ref{sec_dantzig_wolfe_decomposition_subsec_warmstarting_the_column_generation_procedure}.
Given a test instance, we first solve the \gls{LPR} by CPLEX.
Then, the optimal dual values are used as the initial dual values for the column generation procedure.
This method does not require any training.
This is the method used by \citet{Borghettietal2002} and \citet{Schulzeetal2017}.

\paragraph{Column pre-population:}
\revised{
this approach requires training instances to be solved beforehand.
Before the evaluation, for each set of generators, as many training instances as possible are solved to 0.25\% optimality in 24 hours with 8 CPU cores.
The number of training instances solved is reported in Table \ref{tab_number_of_training_instances}.
The training instances are solved using the column generation procedure with the \gls{LPR} initialisation and the local search primal heuristic.
Implementation details are given in Appendix \ref{sec_appendix_implementation_details}.
For each training instance, the problem parameter values and \revisedii{all the generated columns} are saved.
To solve a test instance, its problem parameter values (i.e., the demands) are compared against those of the training instances.
The \revisedii{70} training instances with the closest problem parameter values in terms of the Euclidean distance are selected.
The columns of the selected nearest training instances are retrieved and added to the \gls{RMP} of the test instances.
The \gls{RMP} is then solved without regularisation and the optimal dual values are used as the initial dual values.
\revisediii{In our preliminary experiments, we observed that including the pre-populated columns in the \gls{RMP} increased the \gls{RMP} solution times significantly and degraded the overall performance.
We therefore discard these columns.}
The number of nearest neighbours used, i.e., \revisedii{70}, was chosen by exhaustive grid search on a set of validation \gls{UC} instances which were different from the test instances and the training instances.
That is, the performance of the column pre-population initialisation with different numbers of neighbours was evaluated and the number of neighbours with the best performance was chosen.
}

\begin{table}
\TABLE
{Number of training instances solved\label{tab_number_of_training_instances}}
{%
\begin{tabular}{rr}
\toprule
instance size  & number of training instances solved \\
\midrule
  200  &  25,660 \\
  600  &  12,023 \\
 1000  &   8,588 \\
\bottomrule
\end{tabular}
}{}
\end{table}

\subsubsection{Machine Learning Methods based on Single-Sampling Training.}

The next three methods are based on the single-sampling training and require a dataset of the optimal dual values of training instances.
In our implementation, the datasets created for the column pre-population initialisation (where the optimal dual values are also stored) are used.

\paragraph{Nearest Neighbour:}
given a test instance, the four training instances with the closest problem parameter values are chosen, and the mean of the corresponding optimal dual values are used as initial dual values.
\revised{%
The number of nearest neighbours (i.e., four) is chosen by grid search on validation \gls{UC} instances, in the same way as the column pre-population initialisation.
}

\paragraph{Random forest:}
\revised{%
a random forest model (\citet{Briantetal2008}) is trained using the dataset to predict the optimal dual values given problem parameter values $\omega$.
The training of the models took less than 30 seconds in all cases.
We used the default hyperparameter values of scikit-learn.
}

\paragraph{Neural Network (single-sampling):}
\revised{
the neural network model consists of 4 hidden layers of 1000 units per layer, with skip connections between hidden layers and tanh as \revisedii{an} activation function.
The hyperparameters are chosen based on the performance on validation instances.
For more details on model architecture and training procedure, see Appendix \ref{sec_appendix_hyperparameters_on_neural_networks}.
The time to train the neural network models was less than 15 minutes in every case.
}

\subsubsection{Machine Learning Methods based on Double-Sampling Training.}

The final method is based on the double-sampling training.

\paragraph{Neural Network (double-sampling):}
the initialisation method based on a neural network model with the double-sampling training is implemented as described in Section \ref{subsec_machine_learning_model_based_on_double_sampling}.
The neural network model has the same architecture as the one used for the single-sampling training.
For each set of generators, a single neural network model is trained with the same training time (24 hours using 8 CPU cores).
See Appendix \ref{sec_appendix_hyperparameters_on_neural_networks} for more detail on the training procedure.

\subsection{Evaluation}
\label{sec_neumerical_experiments_subsec_evaluation}

\revised{%
The initialisation methods based on machine learning models are trained to predict dual values which give tight dual lower bounds.
To evaluate the performance of the methods, we first evaluate the dual lower bound at the dual values output by the initialisation methods.
Then, in the next experiment, the column generation procedure is run using the initialisation method but without any primal heuristics and the time to find a near-optimal dual lower bound is measured.
In the final experiment, the column generation procedure is run with primal heuristics and the time to find a near-optimal primal feasible solution is measured.
}

\subsubsection{\revised{Evaluation: Dual Lower Bound}}
\label{sec_neumerical_experiments_subsec_evaluation_dual_lower_bound}

\revised{%
To run the evaluations, 100 test instances of each size were created.
\revisedii{For each test instance, CPLEX was used for two hours to solve it and} the best lower bound found within the time limit was saved (we note that when there is an integrality gap between the Lagrangian lower bound and the optimal objective value if given sufficient time CPLEX will find a tighter lower bound than the column generation procedure).
Then, the dual lower bounds computed at the dual values output by the initialisation methods were calculated by \eqref{problem_eq_dw_lower_bound} and compared against the lower bound found by CPLEX.
}
Table \ref{tab_lower_bound} shows the tightness of the dual lower bounds and the required time (all times in this section are wall clock times) to run the initialisation methods.
The tightness of the dual lower bounds is reported as the average percentage gaps between the bounds.
For comparison, the average gap for the objective value of the \gls{LPR} (which is also a valid lower bound) is shown in the table as well (labelled as \gls{LPR} objective value).

Clearly, coldstart yields poor lower bounds.
\revisedii{The performance of \gls{LPR} initialisation is significantly better than coldstart.
The column pre-population initialisation method gives even better lower bounds, and in particular, is the best on 200-generator instances.
However, the performance is poorer on the largest test instances.
The lower bounds of the methods based on machine learning are comparable, and they give the best lower bounds on large test instances.
Among these methods, the neural network with the double-sampling training performs best on the two largest instances.
}
\revised{%
The computational time of the column pre-population and \gls{LPR} initialisation grows significantly as the problem size increases.
}
On the other hand, the computational time of the other methods remains small.
We note that \revised{as discussed in Section \ref{sec_dantzig_wolfe_decomposition_subsec_initial_dual_lower_bounds_in_cg}} the optimal \gls{LPR} objective value gives a lower bound, but we see from Table \ref{tab_lower_bound} that the dual lower bound evaluated at the optimal dual yields significantly tighter lower bounds.

\begin{table}
\TABLE
{Tightness of initial lower bound $\text{lb}_1$ (\%) and computational time (seconds). The tightness of initial lower bounds $\text{lb}_1$ is measured by comparing it with the best known lower bound $\text{lb}_{\text{CPLEX}}$ obtained by running CPLEX for 2 hours and is reported as the percentage of $\text{lb}_{\text{CPLEX}}$. The values reported in this table are average over 100 test instances.\label{tab_lower_bound}}
{%
\begin{tabular*}{15cm}{@{\extracolsep{\fill}} lrrrrrr}
\toprule
& \multicolumn{2}{c}{size: 200} & \multicolumn{2}{c}{600} & \multicolumn{2}{c}{1000} \\
\cline{2-3}
\cline{4-5}
\cline{6-7}
method  &       $\text{lb}_{\text{CPLEX}} - \text{lb}_1$ &   time &       $\text{lb}_{\text{CPLEX}} - \text{lb}_1$ &    time &       $\text{lb}_{\text{CPLEX}} - \text{lb}_1$ &    time \\
\midrule
coldstart                    &  99.823 &      0.0 & 99.875 &       0.0 & 99.882 &      0.0 \\
LPR (dual lower bound)       &   0.133 &      4.6 &  0.121 &      18.1 &  0.092 &     29.6 \\
LPR objective value          &   0.193 &      4.6 &  0.195 &      18.1 &  0.155 &     29.6 \\
column pre-population         & \revisedii{0.019} & \revisedii{15.0} &  \revisedii{0.025} &      \revisedii{24.2} &  \revisedii{0.073} &     \revisedii{34.5} \\
\ghline
nearest neighbour            &   0.041 & $<\!0.1$ &  0.052 &  $<\!0.1$ &  0.059 & $<\!0.1$ \\
random forest                &   0.056 & $<\!0.1$ &  0.062 &  $<\!0.1$ &  0.064 & $<\!0.1$ \\
network (single-sampling)    &   0.047 & $<\!0.1$ &  0.047 &  $<\!0.1$ &  0.048 & $<\!0.1$ \\
network (double-sampling)    &   0.048 & $<\!0.1$ &  0.037 &  $<\!0.1$ &  0.026 & $<\!0.1$ \\
\bottomrule
\end{tabular*}%
}{}
\end{table}

\subsubsection{\revised{Evaluation: Column Generation Procedure without Primal Heuristic}}

\revised{%
Given the dual values computed by the initialisation methods, the column generation procedure is expected to successively tighten the lower bound.
Hence, if the initialisation methods give good dual values, the column generation procedure is more likely to find a near-optimal dual lower bound in a short time.

To verify this point, the column generation procedure was run with the initialisation methods, but without any primal heuristic, until dual values that yield a dual lower bound of a prescribed suboptimality (0.1\%, 0.05\% and 0.025\%) were found or the time limit of 10 minutes was reached.
We note that, without primal heuristics, the column generation procedure does not provide upper bounds, and therefore no measure of suboptimality of the lower bounds is available.
We have monitored the progress of the column generation procedure (i.e., the suboptimality of the computed lower bounds) by comparing with the best lower bound computed by CPLEX beforehand as explained in Section \ref{sec_neumerical_experiments_subsec_evaluation_dual_lower_bound}.

Table \ref{tab_performance_without_ph} shows the average computational time and the average number of iterations to find the near-optimal dual lower bound or until the time limit of 10 minutes was reached, and the number of problems for which the column generation procedure found a near-optimal dual lower bound within the time limit.
The computational time reported in this table includes the time to run the initialisation methods such as solving \gls{LPR} as well as the time to solve the \gls{RMP} and the pricing subproblems.
For the instances that are not solved within the time limit, the time is set to be 10 minutes and the number of iterations reached by the 10-minute time limit is used.

\revisedii{%
Clearly, coldstart initialisation methods are the slowest, followed by the \gls{LPR} initialisation method.
The methods based on machine learning are the best.
}
On 600-generator instances, the neural network with the double-sampling training is one of the best methods, and on 1000-generator instances, it gives the best performance on average (the standard error of the computational time was less than 4 seconds).
\revisediii{%
There is a correlation between the results in Table 2 and Table 3: if
an initialisation method outputs dual values with a tight lower bound,
the column generation procedure tends to require fewer iterations.
For
example in the 200-generator and 600-generator instances, the column
pre-population initialisation outputs dual values with the tightest
lower bound and the column generation procedure requires a smaller
number of iterations than the other methods.
However, finding the initial dual values using the column pre-population approach is very slow compared to double-sampling initialisation and as a result, the solution obtained using the double-sampling initialisation is on average significantly faster.
}
}

\begin{table}
\TABLE
{Performance of the column generation procedure without primal heuristics. The three columns labelled as ``time'', ``iters'' and ``solved'' show average computational time (seconds) and iterations, and the number of problems solved within the time limit (problems where the column generation procedure found a near-optimal lower bound), respectively.\label{tab_performance_without_ph}}
{%
\setlength\tabcolsep{2pt}  
\begin{tabular*}{15cm}{@{\extracolsep{\fill}} rr rrr rrr rrr}
  \toprule
   &
   &
   \multicolumn{3}{c}{0.1\% optimality} &
   \multicolumn{3}{c}{0.05\% optimality} &
   \multicolumn{3}{c}{0.025\% optimality} \\
  \cline{3-5}
  \cline{6-8}
  \cline{9-11}
  size &
  method &
  \multicolumn{1}{c}{time} & \multicolumn{1}{c}{iters} & solved &
  \multicolumn{1}{c}{time} & \multicolumn{1}{c}{iters} & solved &
  \multicolumn{1}{c}{time} & \multicolumn{1}{c}{iters} & solved \\
\midrule
200 & coldstart                       &   95.9 &        18.2 &      100 &  131.0 &        21.5 &       99 &  181.6 &        25.6 &       97 \\
    & LPR                             &   15.1 &         3.2 &      100 &   23.5 &         5.5 &      100 &   35.8 &         8.4 &      100 \\
    & column pre-population           &  \revisedii{19.9} & \revisedii{1.0} & \revisedii{100} & \revisedii{19.9} & \revisedii{1.0} & \revisedii{100} & \revisedii{21.2} & \revisedii{1.4} & \revisedii{100} \\
\cghline
    & nearest neighbour               &    4.3 &         1.0 &      100 &    8.0 &         2.0 &      100 &   21.6 &         5.3 &      100 \\
    & random forest                   &    4.3 &         1.1 &      100 &    9.8 &         2.5 &      100 &   27.6 &         7.0 &      100 \\
    & network (single-sampling)       &    4.1 &         1.0 &      100 &    6.9 &         1.8 &      100 &   17.7 &         4.4 &      100 \\
    & network (double-sampling)       &    4.7 &         1.1 &      100 &    9.3 &         2.2 &      100 &   25.3 &         6.2 &      100 \\
\midrule
600 & coldstart                       &  422.3 &        15.6 &       76 &  494.2 &        17.3 &       52 &  540.1 &        18.3 &       40 \\
    & LPR                             &   44.1 &         2.8 &      100 &   62.2 &         4.5 &      100 &   78.2 &         5.9 &      100 \\
& column pre-population & \revisedii{37.9} & \revisedii{1.1} & \revisedii{100} & \revisedii{39.1} & \revisedii{1.2} & \revisedii{100} & \revisedii{42.8} & \revisedii{1.5} & \revisedii{100} \\
\cghline
    & nearest neighbour               &   12.7 &         1.2 &      100 &   24.5 &         2.4 &      100 &   50.7 &         4.7 &      100 \\
    & random forest                   &   12.3 &         1.2 &      100 &   29.3 &         2.9 &      100 &   55.8 &         5.3 &      100 \\
    & network (single-sampling)       &   11.4 &         1.1 &      100 &   19.8 &         1.9 &      100 &   39.1 &         3.7 &      100 \\
    & network (double-sampling)       &   12.7 &         1.1 &      100 &   21.1 &         1.8 &      100 &   37.0 &         3.1 &      100 \\
\midrule
1000 & coldstart                      &  450.3 &        13.3 &       84 &  507.4 &        14.9 &       73 &  544.4 &        16.1 &       59 \\
     & LPR                            &   64.2 &         2.1 &      100 &   95.1 &         3.8 &      100 &  117.4 &         5.0 &      100 \\
& column pre-population & \revisedii{60.6} & \revisedii{1.4} & \revisedii{100} & \revisedii{67.7} & \revisedii{1.8} & \revisedii{100} & \revisedii{76.1} & \revisedii{2.2} & \revisedii{100} \\
\cghline
     & nearest neighbour              &   23.8 &         1.5 &      100 &   45.4 &         2.7 &      100 &   76.8 &         4.4 &      100 \\
     & random forest                  &   21.4 &         1.3 &      100 &   50.1 &         3.0 &      100 &   83.8 &         4.8 &      100 \\
     & network (single-sampling)      &   19.2 &         1.2 &      100 &   35.5 &         2.2 &      100 &   62.4 &         3.7 &      100 \\
     & network (double-sampling)      &   18.6 &         1.0 &      100 &   26.4 &         1.4 &      100 &   48.1 &         2.5 &      100 \\
\bottomrule
  \end{tabular*}
}{}
\end{table}

\subsubsection{\revised{Evaluation: Column Generation Procedure with Primal Heuristic}}

\revised{%
In the previous experiments, the initialisation methods were evaluated in terms of the quality of the dual lower bounds.
This is a natural criterion to evaluate the methods since they are trained to maximise the dual lower bound.
However, in practice, the primary interest is to solve the original \gls{UC} instance by finding a good feasible solution.
To see the performance of the initialisation methods in a more practical setup the column generation procedure was run with primal heuristics until a solution within 1\%, 0.5\% and 0.25\% suboptimality was found or the time limit of 10 minutes was reached.
The description of the primal heuristic used is given in Appendix \ref{sec_appendix_implementation_details}.
}

Table \ref{tab_performance} shows the average computational time and the average number of iterations to close the optimality gap or to reach the time limit of 10 minutes, and the number of problems solved within the time limit.
For comparison, we also solve the extensive form \eqref{problem_eq_mip} to the same tolerances without decomposition (i.e., by branch and bound) using CPLEX.

\begin{table}
\TABLE
{Performance of the column generation procedure with primal heuristics. The three columns labelled as ``time'', ``iters'' and ``solved'' show average computational time (seconds) and iterations, and the number of problems solved within the time limit (problems where the column generation procedure found a near-optimal primal solution), respectively.\label{tab_performance}}
{%
\setlength\tabcolsep{2pt}  
\begin{tabular*}{15cm}{@{\extracolsep{\fill}} rr rrr rrr rrr}
  \toprule
   &
   &
   \multicolumn{3}{c}{1\% optimality} &
   \multicolumn{3}{c}{0.5\% optimality} &
   \multicolumn{3}{c}{0.25\% optimality} \\
  \cline{3-5}
  \cline{6-8}
  \cline{9-11}
  size &
  method &
  \multicolumn{1}{c}{time} & \multicolumn{1}{c}{iters} & solved &
  \multicolumn{1}{c}{time} & \multicolumn{1}{c}{iters} & solved &
  \multicolumn{1}{c}{time} & \multicolumn{1}{c}{iters} & solved \\
  \midrule
200  & CPLEX                        &  261.4 &           - &       98 &  271.3 &           - &       98 &  392.6 &           - &       94 \\
\cghline
     & coldstart                    &  106.7 &        16.4 &      100 &  151.3 &        20.4 &      100 &  237.6 &        27.2 &       99 \\
     & LPR                          &    8.3 &         1.1 &      100 &   18.8 &         3.4 &      100 &   60.8 &        11.6 &      100 \\
& column pre-population & \revisedii{20.0} & \revisedii{1.0} & \revisedii{100} & \revisedii{21.2} & \revisedii{1.3} & \revisedii{100} & \revisedii{51.8} & \revisedii{7.0} & \revisedii{100} \\
\cghline
     & nearest neighbour            &    4.8 &         1.0 &      100 &    8.3 &         1.7 &      100 &   34.2 &         6.8 &      100 \\
     & random forest                &    4.9 &         1.1 &      100 &    9.2 &         2.0 &      100 &   44.4 &         8.8 &      100 \\
     & network (single-sampling)    &    5.4 &         1.1 &      100 &   12.5 &         2.6 &      100 &   40.9 &         8.1 &      100 \\
     & network (double-sampling)    &    4.8 &         1.0 &      100 &    8.4 &         1.7 &      100 &   36.6 &         7.3 &      100 \\
  \hline
600  & CPLEX                        &  593.2 &           - &        6 &  594.2 &           - &        5 &  597.9 &           - &        4 \\
\cghline
     & coldstart                    &  416.2 &        13.9 &       83 &  505.3 &        15.9 &       53 &  564.0 &        17.4 &       25 \\
     & LPR                          &   29.1 &         1.1 &      100 &   44.9 &         2.2 &      100 &   93.5 &         5.7 &      100 \\
& column pre-population & \revisedii{39.0} & \revisedii{1.1} & \revisedii{100} & \revisedii{42.0} & \revisedii{1.3} & \revisedii{100} & \revisedii{59.0} & \revisedii{2.4} & \revisedii{100} \\ 
\cghline
     & nearest neighbour            &   12.7 &         1.0 &      100 &   21.8 &         1.6 &      100 &   51.5 &         3.9 &      100 \\
     & random forest                &   13.1 &         1.0 &      100 &   22.8 &         1.8 &      100 &   56.5 &         4.3 &      100 \\
     & network (single-sampling)    &   14.3 &         1.1 &      100 &   28.8 &         2.1 &      100 &   53.0 &         3.9 &      100 \\
     & network (double-sampling)    &   13.3 &         1.0 &      100 &   16.8 &         1.2 &      100 &   44.5 &         3.1 &      100 \\
  \hline
1000 & CPLEX                        &  600.0 &           - &        0 &  600.0 &           - &        0 &  600.0 &           - &        0 \\
\cghline
     & coldstart                    &  508.9 &        11.8 &       63 &  559.5 &        13.1 &       36 &  590.0 &        14.0 &       15 \\
     & LPR                          &   46.5 &         1.0 &      100 &   72.2 &         2.0 &      100 &  131.2 &         4.6 &      100 \\
& column pre-population & \revisedii{61.9} & \revisedii{1.3} & \revisedii{100} & \revisedii{74.3} & \revisedii{1.8} & \revisedii{100} & \revisedii{112.6} & \revisedii{3.4} & \revisedii{100} \\
\cghline
     & nearest neighbour            &   19.5 &         1.0 &      100 &   43.3 &         2.0 &      100 &   90.7 &         4.2 &      100 \\
     & random forest                &   20.8 &         1.1 &      100 &   43.5 &         2.0 &      100 &   98.0 &         4.4 &      100 \\
     & network (single-sampling)    &   18.9 &         1.0 &      100 &   49.5 &         2.3 &      100 &   94.0 &         4.4 &      100 \\
     & network (double-sampling)    &   20.1 &         1.0 &      100 &   24.9 &         1.2 &      100 &   65.9 &         2.9 &      100 \\
    \bottomrule
  \end{tabular*}
}{}
\end{table}

These results show that solving the problems without decomposition is the slowest.
\revisedii{%
The other methods show overall the same trend as Table \ref{tab_performance_without_ph}.
Solving the problem with coldstart is by far the slowest in every case.
The column pre-population initialisation is slower than the remaining methods with the loose tolerance (1.0\%), but faster than the \gls{LPR} initialisation method with the tighter tolerance (0.5\% or 0.25\%).
}%
The methods based on machine learning are faster in all cases than the \gls{LPR}, \revised{column pre-population and coldstart initialisation methods}.
\revised{%
When the problem is small or the tolerance is loose, the neural network with the double-sampling training is one of the best methods.
On the large instances (600 or 1000 generators) with the tight tolerance (0.5\% or 0.25\%) it gives the best performance.
}
\revisediii{%
We note that, as we saw in Table 3, there is again a correlation between the tightness of initial lower bounds computed by an initialisation method and the number of iterations required by the column generation procedure.
For example, comparing the pre-populate
and the double sample results only 2 cases of 9 (the 1\% and 0.5\% for 600 generators) are against the trend and in these the deviation is minor.
}

To observe the effect of the quality of the initial dual values on the time taken by the column generation procedure more clearly, in Figure \ref{fig_time_vs_lb} the total computational time (0.25\% tolerance) for each 1000-generator test instance is plotted against the tightness of the initial lower bound.
The initialisation methods based on the \gls{LPR} and the neural network model with the double-sampling training are compared in the plot.
We observe a correlation between the two metrics.
That is, if the lower bound computed in the first iteration of the column generation procedure is tighter then the total computational time required by the column generation procedure with primal heuristic tends to be smaller.

\begin{figure}
    \begin{center}  
    \small
\begin{tikzpicture}
\begin{axis}[
    xlabel={Total solution time (seconds)},
    ylabel={Tightness of initial lower bounds (\%)},
    xmin=0,
    xmax=300,
    ymin=0,
    ymax=0.21,
    ytick={0.00, 0.05, 0.10, 0.15, 0.20},
    yticklabels={0.00, 0.05, 0.10, 0.15, 0.20},
    legend pos=south east,
    legend cell align=left,
    legend style={
      /tikz/every even column/.append style={column sep=5mm}
    },
    width=.65\textwidth,
    height=0.30\textheight,
]


\tikzset{
  style_online_network/.style={mark=o,mark size=1.4pt,mark options={draw=red,fill=red}},
  style_lpr/.style={mark=x,mark size=2.4pt,mark options={draw=black,fill=black}},
};
\tikzset{
  lr_style_online_network/.style={red,very thin},
  lr_style_lpr/.style={black,very thin},
};
\pgfplotsinvokeforeach{online_network,lpr} {
    \edef\temp{\noexpand
  \addplot+ [
        style_#1,
        only marks,
  ]
      table [
        x=total_time,
        y=first_lb,
        col sep=comma
      ] {./data/warmstart_plot_data_100/1000_#1.csv};
  \noexpand
  \addplot+ [
      lr_style_#1,
      no marks,
      forget plot
  ]
      table [
        x=total_time,
        y={create col/linear regression={y=first_lb, x=total_time}},
        col sep=comma
      ] {./data/warmstart_plot_data_100/1000_#1.csv};
    };  
    \temp
}
\legend{Network (double-sampling),LPR};
\end{axis}
\end{tikzpicture}
\end{center}  
\caption{Tightness of the initial lower bound vs the total computational time required by the column generation procedure on 1000-generator test instances}\label{fig_time_vs_lb}
\end{figure}

A method to generate initial dual values must balance the time to train, the quality of the dual values and the computational time.
One extreme example is the \gls{LPR} initialisation.
It does not require any offline training.
However, it produces dual values with a relatively loose lower bound and the solution time grows as the problem size increases.
The methods based on machine learning models need offline training, but they run very quickly when solving a new instance and output dual values with a tight lower bound, resulting in a significant reduction of the time required to solve new instances.

To clarify this point, Table \ref{tab_solution_time_breakdown} shows a breakdown of the average computational time of the column generation procedures applied to the 1000-generator instances for the 100 test cases.
The column labelled `Initialisation' shows the time required to run the initialisation methods.
For the \gls{LPR} initialisation method this is the time to solve the \gls{LPR}, while for the column pre-population initialisation method, this is the time to solve the unregularised \gls{RMP}.
For the methods based on machine learning, this is the time to evaluate the models.
In the case of $1\%$ tolerance, the time required to solve the \gls{LPR} is longer than the sum of the time spent on the other routines.
However, as the tolerance becomes tighter, the number of iterations and the time spent in the other routines increase, and by $0.25\%$ tolerance, the \gls{LPR} initialisation time is only 25\% of the total time.
\revisedii{We see a similar trend if we use the column pre-population initialisation.}

\begin{table}
\TABLE
{Breakdown of the average computational time (seconds) for 1000-generator case.\label{tab_solution_time_breakdown}}
{
  \begin{tabular}{lr rrrr}
\toprule
tolerance      &          & Initialisation & RMP &  Subproblem & Primal Heuristic \\
\midrule
1\%    & coldstart                   &     - &  238.9 &  186.9 &  90.0 \\
       & LPR                         &  29.6 &    0.1 &   12.2 &   4.7 \\
       & column pre-population       &   \revisedii{34.5} & \revisedii{0.7} & \revisedii{22.2} & \revisedii{4.5} \\ 
\cdashline{2-6}
       & nearest neighbour           &   0.0 &    0.1 &   15.0 &   4.4 \\
       & random forest               &   0.0 &    0.1 &   15.6 &   5.0 \\
       & network (single-sampling)   &   0.0 &    0.1 &   15.2 &   3.6 \\
       & network (double-sampling)   &   0.0 &    0.1 &   16.6 &   3.4 \\
\hline
0.5\%  & coldstart                   &     - &  269.2 &  209.7 &  93.7 \\
       & LPR                         &  29.6 &    1.0 &   28.9 &  12.8 \\
       & column pre-population       &   \revisedii{34.5} & \revisedii{2.4} & \revisedii{30.6} & \revisedii{6.7} \\                                                                                                 
\cdashline{2-6}
       & nearest neighbour           &   0.0 &    0.7 &   31.1 &  11.5 \\
       & random forest               &   0.0 &    0.9 &   30.5 &  12.1 \\
       & network (single-sampling)   &   0.0 &    0.9 &   35.4 &  13.2 \\
       & network (double-sampling)   &   0.0 &    0.2 &   20.0 &   4.8 \\
\hline
0.25\% & coldstart                   &     - &  289.9 &  222.8 &  95.9 \\
       & LPR                         &  29.6 &    4.8 &   69.0 &  27.8 \\
       & column pre-population       &   \revisedii{34.5} & \revisedii{5.5} & \revisedii{59.5} & \revisedii{13.0} \\
\cdashline{2-6}
       & nearest neighbour           &   0.0 &    2.8 &   65.1 &  22.8 \\
       & random forest               &   0.0 &    3.5 &   68.3 &  26.3 \\
       & network (single-sampling)   &   0.0 &    2.9 &   67.7 &  23.4 \\
       & network (double-sampling)   &   0.0 &    1.2 &   50.0 &  14.6 \\
\bottomrule
  \end{tabular}
}{}
\end{table}

\section{Conclusion}
\label{sec_conclusion}

We have investigated the use of machine learning techniques to accelerate Dantzig-Wolfe decomposition with a column generation procedure to solve parametrized \gls{UC} problems.
\revised{%
In particular, we have proposed the use of machine learning models to compute the initial dual values for the column generation procedure.
We have trained machine learning models so that they output dual values that yield tight dual lower bounds.
In contrast to previous approaches which construct primal solutions directly, our approach can deliver an accurate solution as well as obtain its suboptimality bound, from the column generation procedure.
}

\revised{%
We have considered two approaches to training a machine learning model.
The first approach, the single-sampling training, was to directly predict the optimal dual values given the problem parameter values as a standard regression problem.
This requires a large dataset of the optimal dual values of training instances.
The second approach, the double-sampling training, was based on the decomposable structure of the dual function.
}

We have first evaluated the performance of the machine learning models by computing the dual lower bounds at the dual values generated by the models.
As benchmarks, we have used the dual lower bounds using the optimal dual values to the \gls{LPR}, the optimal dual values to the \gls{RMP} with pre-populated columns and coldstart with dual values $y = 0$.
The coldstart initialisation always has the worst dual lower bounds, while the \gls{LPR} generates dual values with relatively tight dual lower bounds.
\revisedii{%
The column pre-population initialisation gives the best dual lower bounds on 200-generator instances, but the performance is poorer on the largest test instances.
}
The methods based on machine learning yield dual values with significantly tighter dual lower bounds, compared to the bound computed by the \gls{LPR}, \revisedii{and the best bounds on the larger test instances}.
The computational time to solve the \gls{LPR} or the \gls{RMP} with pre-populated columns grows significantly as the instance size becomes larger.
However, the computational time to evaluate the machine learning models remains negligible even with large instances.

\revised{%
To see the effectiveness of the above initialisation methods to warmstart the column generation procedure, we have evaluated the performance of the column generation procedure when warmstarted with the dual values obtained by each initialisation method.
For this reason, the column generation procedure was run without any primal heuristics.
In this evaluation, the time required by the column generation procedure to find near-optimal dual lower bounds was measured.
The results revealed that the column generation procedure initialised with any of all the machine learning based initialisation methods can find a near-optimal dual lower bound more quickly than with the \gls{LPR}, column pre-population or coldstart initialisation.
}

\revised{%
We have also evaluated the performance of the warmstarted column generation procedure in a more practical setup where a near-optimal primal solution (generator schedule) is sought.
In this experiment, to find primal feasible solutions, the column generation procedure was run with primal heuristics and we observed the time required to find a primal solution of a prescribed suboptimality.
}
In the numerical experiments, we observed that solving the \gls{UC} problem with decomposition was always faster than solving the problem without decomposition using CPLEX.
We further noted that warmstarting the column generation procedure successfully reduced the number of iterations and overall computational time to find a solution of prescribed suboptimality and this was especially dramatic for the initialisation methods based on machine learning.
We observed that the initialisation methods that generate tighter initial Lagrangian lower bounds produce a better performance of the column generation procedure with primal heuristics.
For example, the initialisation methods based on machine learning outperformed the \gls{LPR} for all problem sizes.
The numerical experiments also showed that the methods based on machine learning scale well and can be effective for solving large-scale \gls{UC} problems.
\revised{%
In particular, the neural network initialisation with the double-sampling training was usually the best, especially on the large instances (Table \ref{tab_lower_bound}).
One possible reason for the strong performance of the neural network with the double-sampling training compared to the other machine learning models is the efficiency of the training.
The double-sampling training uses a stochastic gradient method which exploits the structure of the problem.
With the training time we used (24 hours), we do not observe a significant difference between the methods on 200-generator instances, however, on larger instances, there were differences between them.
}

An interesting area for further study is the \gls{UC} problem with stochastic demand forecasts.
In many solution methods, \gls{UC} problems with deterministic demand forecasts appear as subproblems and are solved repeatedly.
See \citet{Takritietal1996} and \citet{Schulzeetal2017} for instance.
Typically the subproblems share the same problem structure and only differ in the cost coefficients and the demand forecasts.
This paper considers only the case where the perturbations are on the demand but not on the cost.
However, the technique we developed is extendable to the stochastic case and may be equally effective.

\ACKNOWLEDGMENT{%
The authors are deeply grateful to two anonymous referees: their valuable and insightful discussions led to a significant improvement of the paper.
}

%
%
%
\begin{APPENDICES}

\section{Problem Formulation}
\label{sec_appendix_problem_formulation}

We closely follow one of the standard formulations in literature, referred to as the 3-binary variable formulation by \citet{OstrowskiEtAl2012}, and formulate the following constraints:
\begin{itemize}
\item \textbf{Load balance}: Generators have to meet all the demand in each time period (generation shedding at 0 cost is allowed).
\item \textbf{Reserve}: To deal with contingencies, it is required to keep a sufficient amount of backup in each time period, which can be activated quickly.
\item \textbf{Power output bounds}: Each generator's power output has to be within its limit.
\item \textbf{Ramp rate bounds}: Generators can only change their outputs within the ramp rates.
\item \textbf{Minimum up/downtime}: If switched on (off), each generator has to stay on (off) for a given minimum period.
\end{itemize}

The formulation of the model is as follows.

\begin{itemize}
\item{Parameters}
\begin{itemize}
\item $G$: number of generators
\item $T$: number of time periods where decisions are taken
\item $C^{\mathrm{nl}}_{g}$: no-load cost of generator $g$
\item $C^{\mathrm{mr}}_{g}$: marginal cost of generator $g$
\item $C^{\mathrm{up}}_{g}$: startup cost of generator $g$
\item $P^{\max/\min}_{g}$: maximum/minimum generation limit of generator $g$
\item $P^{\mathrm{ru}/\mathrm{rd}}_{g}$: operating ramp up/down limits of generator $g$
\item $P^{\mathrm{su}/\mathrm{sd}}_{g}$: startup/shutdown ramp limits of generator $g$
\item $T^{\mathrm{u}/\mathrm{d}}_{g}$: minimum up/downtime of generator $g$
\item $P^{\mathrm{d}}_{t}$: power demand at time $t$
\item $P^{\mathrm{r}}_{t}$: reserve requirement at time $t$
\end{itemize}

\item{Variables}
\begin{itemize}
\item $\alpha_{gt} \in \{0, 1\}$: 1 if generator $g$ is on in period $t$, and 0 otherwise
\item $\gamma_{gt} \in \{0, 1\}$: 1 if generator $g$ starts up in period $t$, and 0 otherwise
\item $\eta_{gt} \in \{0, 1\}$: 1 if generator $g$ shuts down in period $t$, and 0 otherwise
\item $p_{gt} \ge 0$: power output of generator $g$ in period $t$
\end{itemize}

\item Total cost (the objective to be minimised)
\[
\min \sum_{t = 1}^T \sum_{g = 1}^G
\left( C^{\mathrm{nl}}_g \alpha_{gt} + C^{\mathrm{mr}}_g p_{gt} +
C^{\mathrm{up}}_g \gamma_{gt}
\right).
\]

\item Load balance
\begin{equation*}
\sum_{g = 1}^G p_{gt} \ge P^{\mathrm{d}}_{t}
\qquad t = 1, 2, \ldots, T.
\label{eq:uc_first_constraint}
\end{equation*}

\item Reserve
\begin{equation*}
\sum_{g = 1}^G (P^{\max}_{g} \alpha_{gt} - p_{gt})
\ge P^{\mathrm{r}}_t
\qquad t = 1, 2, \ldots, T.
\end{equation*}

\item Power output bounds
\begin{equation*}
P^{\min}_{g} \alpha_{gt} \le p_{gt} \le P^{\max}_{g} \alpha_{gt}
\qquad g = 1, 2, \ldots, G, t = 1, 2, \ldots, T
\end{equation*}

\item Ramp rate bounds
\begin{equation*}
p_{gt} - p_{g \, t-1} \le P^{\mathrm{ru}}_g \alpha_{g \, t-1}
+ P^{\mathrm{su}}_g \gamma_{gt}
\qquad g = 1, 2, \ldots, G, t = 2, 3, \ldots, T.
\end{equation*}
\begin{equation*}
p_{g \, t-1} - p_{gt} \le P^{\mathrm{rd}}_g \alpha_{gt}
+ P^{\mathrm{sd}}_g \eta_{gt}
\qquad g = 1, 2, \ldots, G, t = 2, 3, \ldots, T.
\end{equation*}

\item Minimum up/downtime
\begin{equation*}
\sum_{u=\max\{t-T^\mathrm{u}_g+1, 1\}}^t \gamma_{gu} \le \alpha_{gt}
\qquad g = 1, 2, \ldots, G, t = 1, 2, \ldots, T
\end{equation*}
\begin{equation*}
\sum_{u=\max\{t-T^\mathrm{u}_g+1, 1\}}^t \eta_{gu} \le 1 - \alpha_{gt}
\qquad g = 1, 2, \ldots, G, t = 1, 2, \ldots, T
\end{equation*}

\item Logical constraints (to enforce binaries to work as we expect)
\begin{equation*}
\alpha_{gt} - \alpha_{g \, t-1} = \gamma_{gt} - \eta_{gt}
\qquad g = 1, 2, \ldots, G, t = 2, 3, \ldots, T
\end{equation*}
\begin{equation*}
1 \ge \gamma_{gt} + \eta_{gt}
\qquad g = 1, 2, \ldots, G, t = 1, 2, \ldots, T
\label{eq:uc_last_constraint}
\end{equation*}
\end{itemize}

\section{Implementation Details}
\label{sec_appendix_implementation_details}

In the numerical experiments, a regularised column generation procedure is used.
Initially, the regularisation centre is set to the dual values given by the initialisation method.
In each iteration, a lower bound is evaluated using the current dual values and the regularisation centre is updated to the current dual values if the lower bound gets improved.
Furthermore, if the lower bound improves, the regularisation parameter is divided by two, and otherwise multiplied by two.
This closely follows the implementation of the column generation procedure described by \citet{Schulzeetal2017}.

In every iteration, after the pricing subproblems are solved, a primal heuristic based on local search is run.
Given the solutions to the pricing subproblems, this primal heuristic checks the feasibility and if necessary switches on the cheapest available generators to make the solution feasible.
Note that the solutions to the pricing subproblems are feasible generator schedules and infeasibility only arises from an insufficient generation capacity to meet the demand or insufficient reserve.
This primal heuristic loosely follows that proposed by \citet{Guanetal1992}.
In our experiments, we observe that it usually finds near-optimal primal solutions when the dual values get close to optimal.
If the column generation procedure using the local search primal heuristic fails to find a primal feasible solution satisfying the given optimality tolerance within 30 iterations, we use in addition a column combination primal heuristic.
Then in the subsequent iterations, we use both of the primal heuristics (the two primal heuristics are independent of each other).
The column combination primal heuristic is a popular primal heuristic in the column generation procedure in general.
The idea of this heuristics is to solve (\ref{problem_eq_mip}) with restricted patterns of solutions, referred to as restricted master IP by \citet{Vanderbeck2005}.
In the $k$th iteration, we replace $X_s$ in (\ref{problem_eq_mip}) with $\{x_s^{(k - l)}\}_{l = 0, 1, 2}$
where $\{x_s^{(k - l)}\}_{l = 0, 1, 2}$ are solutions to the $s$th pricing subproblem found in the current iteration and or the previous two iterations.
The resulting problem is still a mixed-integer programme but the solution space is much smaller than the original problem.
In our numerical experiments, we observed that the problem is feasible without the need to add artificial variables (columns).

\section{Hyperparameters for the Neural Networks}
\label{sec_appendix_hyperparameters_on_neural_networks}

The model consists of 4 hidden layers of 1000 units per layer, with skip connections between hidden layers as described by \citet{DeAndSmith2020}.
The structure of the network is visualised in Figure \ref{sec_appendix_hyperparameters_fig_network}.
The tanh function is used as an activation function.
\revisedii{These hyperparameters were chosen by grid search, as discussed in the end of this section.}
The weights in the linear transformation are initialised based on the methods of \citet{GlorotAndBengio2010}.
All biases and the weights on the residual connections are initialised to zero.

\begin{figure}[hbtp]
\centering  
\begin{tikzpicture}
\node (input) at (0, 0) {};
\node [rotate=-90] (inputlabel) at (input) {input};

\node (lin1) [right=of input] {};
\node[shape=rectangle, draw, rotate=-90] (linlabel1) at (lin1) {linear};
\node (tanh1) [node distance=0cm and 0.5cm, right=of lin1] {};
\node[shape=rectangle, draw, rotate=-90] (tanhlabel1) at (tanh1) {activation};

\node (lin2) [right=of tanh1] {};
\node[shape=rectangle, draw, rotate=-90] (linlabel2) at (lin2) {linear};
\node (tanh2) [node distance=0cm and 0.5cm, right=of lin2] {};
\node[shape=rectangle, draw, rotate=-90] (tanhlabel2) at (tanh2) {activation};
\node (alpha2) [node distance=0cm and 0.5cm, right=of tanh2] {};
\node[shape=rectangle, draw, rotate=-90] (alphalabel2) at (alpha2) {$\alpha$};

\node (lin3) [right=of alpha2] {};
\node[shape=rectangle, draw, rotate=-90] (linlabel3) at (lin3) {linear};
\node (tanh3) [node distance=0cm and 0.5cm, right=of lin3] {};
\node[shape=rectangle, draw, rotate=-90] (tanhlabel3) at (tanh3) {activation};
\node (alpha3) [node distance=0cm and 0.5cm, right=of tanh3] {};
\node[shape=rectangle, draw, rotate=-90] (alphalabel3) at (alpha3) {$\alpha$};

\node (tanh4) [right=of alpha3] {};
\node[shape=rectangle, draw, rotate=-90] (tanhlabel4) at (tanh4) {linear};

\node (output) [right=of tanh4] {};
\node [rotate=-90] (outputlabel) at (output) {output};

\draw[->] (inputlabel) -- (linlabel1);
\draw[->] (linlabel1) -- (tanhlabel1);
\draw[->] (tanhlabel1) -- (linlabel2) node [inner sep=0cm] (firstsecond) [midway] {};
\draw[->] (linlabel2) -- (tanhlabel2);
\draw[->] (tanhlabel2) -- (alphalabel2);
\draw[->] (alphalabel2) -- (linlabel3)
          node [inner sep=0cm] (secondthirdone) [near start] {}
          node [inner sep=0cm] (secondthirdtwo) [near end] {};
\draw[->] (linlabel3) -- (tanhlabel3);
\draw[->] (tanhlabel3) -- (alphalabel3);
\draw[->] (alphalabel3) -- (tanhlabel4) node [inner sep=0cm] (thirdforth) [midway] {};
\draw[->] (tanhlabel4) -- (outputlabel);

\coordinate [below=1.2cm of tanh2] (tanh2below);
\coordinate [below=1.2cm of tanh3] (thirdtanhbelow);

\draw[->] (firstsecond) |- (tanh2below) -| (secondthirdone);
\draw[->] (secondthirdtwo) |- (thirdtanhbelow) -| (thirdforth);

\end{tikzpicture}
\caption{The structure of the neural network model with skip connections. For a more detailed explanation of the architecture, see \citet{DeAndSmith2020}.\label{sec_appendix_hyperparameters_fig_network}}
\end{figure}

\revised{%
For the single-sampling training, we use the dataset of the optimal dual values of training instances.
We split the dataset into a training set (80\%) and a validation set (20\%).
The Adam method (\citet{KingmaAndJimmy2014}) is used to learn the parameters of the neural network.
After each epoch, we evaluate the mean squared error on the validation set.
If it does not improve for successive 4 epochs, we halve the learning rate.
If it does not improve for successive 12 epochs, we terminate the training to avoid overfitting.
}

For the double-sampling training, we also use the Adam method (\citet{KingmaAndJimmy2014}).
The model performance is evaluated every 5 minutes \revisedii{by computing the dual lower bound on validation instances using the output of the neural network.
To this end, 10 validation instances are sampled, which are distinct from the instances used to train the model (i.e. compute \eqref{eq_sampled_gradient_dl_dtheta}) and to test the final performance (reported in Table \ref{tab_performance_without_ph}, \ref{tab_performance}).}
If it fails to improve the performance for successive 15 minutes, the learning rate is divided by 1.5.
\revisedii{
In our experiment, we did not observe overfitting.
However, if this became an issue, we could use regularisation such as dropout, as proposed by \citet{CobbeEtAl2019}.
}

\revisedii{%
As discussed earlier in this section, the number of layers and the activation function were chosen by grid search.
Models with 3, 4 and 5 layers and with tanh and relu were trained on the 200-generator case.
Table \ref{tab_grid_search} reports the performance of the trained models in the same format as Table \ref{tab_lower_bound}.
The instances used in Table \ref{tab_grid_search} are different from the training, validation (used to monitor the progress of the training) and test instances.
\begin{table}
\TABLE
{Tightness of lower bound computed with the output of the neural networks. The tightness is measured by comparing it with the best known lower bound $\text{lb}_{\text{CPLEX}}$ obtained by running CPLEX for 2 hours and is reported as the percentage of $\text{lb}^*$.\label{tab_grid_search}}
{%
\begin{tabular*}{10cm}{@{\extracolsep{\fill}} llr}
\toprule
activation & number of layers &  $\text{lb}_{\text{CPLEX}} - \text{lb}_1$ \\
\midrule
tanh &  3 &   0.0239 \\
     &  4 &   0.0228 \\
     &  5 &   0.0239 \\
\hline
relu &  3 &   0.0462 \\
     &  4 &   0.0480 \\
     &  5 &   0.0479 \\
\bottomrule
\end{tabular*}%
}{}
\end{table}
}

\section{Statistics Relevant to Training}
\label{sec_appendix_statistics_relevant_to_training}

\revisedii{%
Table \ref{tab_difference_of_training_and_test} shows the difference between training and test instances.
For each test instance the commitment decision of the nearest training instance $p^{\text{train}}$ and the commitment decision of the test instance $p^{\text{test}}$ was compared using the Hamming distance:
\[
\text{Difference in terms of primal values (\%)} = \frac{\text{Hamming distance of $p^{\text{train}}$ and $p^{\text{test}}$}}{\text{the number of commitment decisions in $p^{\text{test}}$}} \cdot 100.
\]
The average of the above quantity using 40 test instances was computed and shown in the middle column of Table \ref{tab_difference_of_training_and_test}.
Similarly, the distance between the optimal dual values of the nearest training instances and the test instances was computed.
For each test instance, the nearest training instance is obtained.
Then, the optimal dual values of the training instance and the test instance ($d^{\text{train}}$ and $d^{\text{test}}$, respectively) and the following metric is computed:
\[
\text{Difference in terms of dual values (\%)} = \frac{\|d^{\text{train}} - d^{\text{test}}\|_2}{\| d^{\text{test}} \|_2} \cdot 100.
\]
The average is shown in the right-most column in the table.
}

\begin{table}
\TABLE
{Difference of primal and dual values between training and test instances.\label{tab_difference_of_training_and_test}}
{
\begin{tabular}{lrr}
\toprule
size & difference in terms of primal values (\%) & difference in terms of dual values (\%) \\
\midrule
 200 &  2.93 & 6.85 \\
 600 &  2.88 & 5.99 \\
1000 &  2.80 & 7.09 \\
\bottomrule
\end{tabular}
}{}
\end{table}

\end{APPENDICES}

\ACKNOWLEDGMENT{%
}


\printbibliography


\end{document}